\documentclass{article}
\usepackage{ruble} 
\pdfoutput=1

\def\thetitle{Плоское оригами и длинный рубль}
\def\theauthors{Антон~Петрунин}

\hypersetup{colorlinks=true,
citecolor=black,
linkcolor=black,
anchorcolor=black,
filecolor=black,
menucolor=black,
urlcolor=black,
pdftitle={\thetitle},
pdfauthor={\theauthors}
}

\begin{document}
\title{\thetitle}
\author{\theauthors%
\footnote{К заметке 2008 года добавлен ответ к задаче 5.
\\\hspace*{\parindent}
Большое спасибо Арсению Акопяну, Александру Вострикову, Роберту Лэнгу и Алексею Тарасову за помощь при составлении этой заметки.
}}
\date{}

\maketitle

Задачу, о которой пойдёт речь, придумал 
Владимир Арнольд в 19 лет. 
Сейчас ему уже за 70, он один из самых известных российских математиков.
Известен он много чем, но более всего своими задачами, которых напридумывал уже около тысячи.
Многие из этих задач до сих пор не решены, а чтобы решить некоторые, пришлось развить новые области математики.

Изначально задача формулировалась так: \textit{«можно ли из рубля сложить плоскую фигуру с б\'{о}льшим периметром?»}.
Эта задача была популярна в математическом фольклоре, и в некотором смысле до сих пор не решена.
Мы сформулируем задачу немного по-другому:

\begin{thm}{Задача}
Можно ли сложить квадратный лист бумаги на плоскости так, что периметр полученной фигуры превысит периметр исходного листа?
\end{thm}
\begin{wrapfigure}{r}{55mm}
\noi\begin{lpic}[t(0mm),b(0mm),r(0mm),l(0mm)]{pics/skladka(0.6)}
\lbl{25,28;$M$}
\lbl{75,28;$M'$}
\lbl[r]{23,55;$q$}
\lbl[r]{87,55;$q$}
\end{lpic}
\caption{\label{sgib}Перегиб $M$.}
\end{wrapfigure}

Ответ на эту задачу зависит от того, что понимать под словом «сложить». 
Например, под «складыванием» можно понимать последовательность следующих «перегибаний»: возьмём многоугольник $M$, выберем прямую $q$, секущую $M$, перегнём $M$ вдоль $q$, при этом части $M$ по разные стороны от $q$ «склеятся» в новый многоугольник $M'$ (см. рисунок~\ref{sgib}).

Как следует из упражнения, такими операциями увеличить периметр нельзя. 

\medskip
\noi\textbf{Упражнение.} {\it Докажите, что периметр $M'$ не превосходит периметр $M$. (При решении этого упражнения обратите внимание, что в общем случае прямая $q$ может пересекаться с $M$ по нескольким отрезкам).}
\medskip

При помощи этой операции нельзя отогнуть часть листа, которая лежит на другом. Например, невозможно провести последовательность «перегиба» и «отгиба» как на рисунке~\ref{otgib}.
\begin{figure}
\ \ \ \ \ \ \ \ \begin{lpic}[t(0mm),b(0mm),r(0mm),l(0mm)]{pics/otgib(0.7)}
\end{lpic}
\caption{\label{otgib}Перегиб и отгиб.}
\end{figure} 
Вторая операция увеличивает периметр, хотя он остаётся всё ещё меньше периметра исходного квадрата.
Если под «складыванием»  понимать последовательность таких перегибаний и отгибаний, то ответ на задачу до сих пор никто не знает (хотя, возможно, никто и не хочет узнать).

\section{Японский журавлик}

\begin{figure}[h]
\ \ \ \ \ \ \ \ 
\includegraphics[scale=0.34]{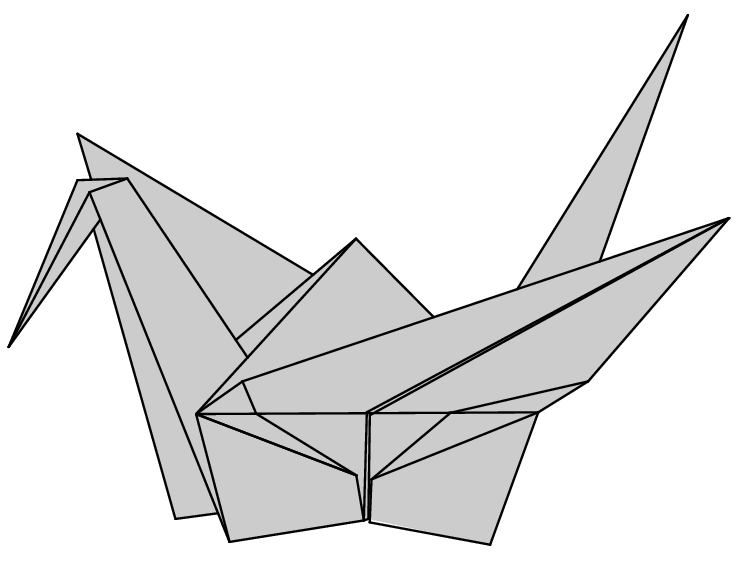}
\hfill
\includegraphics[scale=0.18]{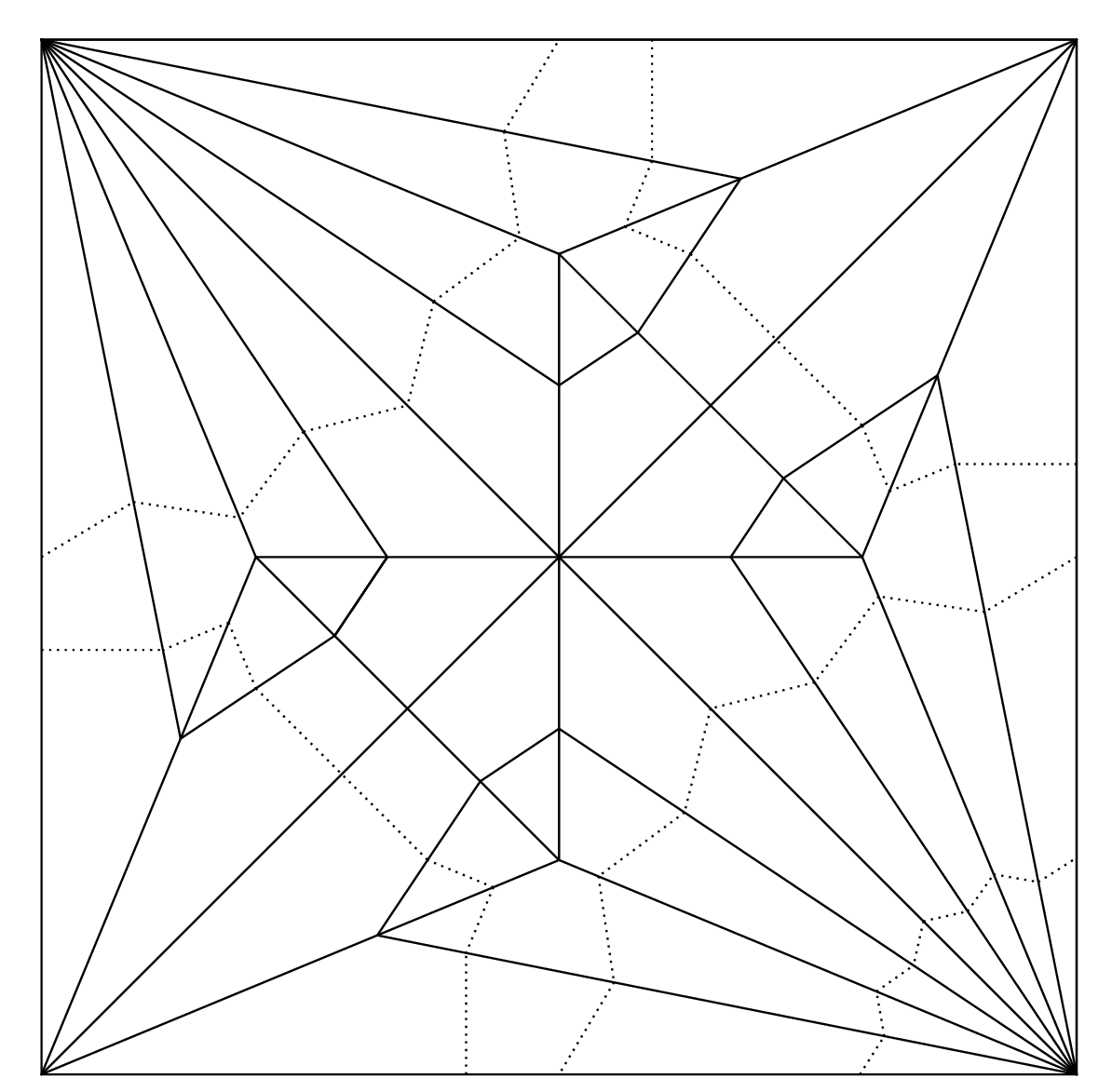}
\ \ \ \ \ \ \ \ 
\caption{Японский журавлик и его развёртка. Точечным пунктиром обозначены складки, которые добавляются при расправлении крыльев, загибании шеи, хвоста и головы (голова --- правый нижний угол, крылья --- левый нижний и правый верхний, хвост --- левый верхний).}
\end{figure}

Под складыванием можно понимать ещё более общую операцию: 
представьте себе, что мы разметили складки заранее и сгибаем лист целиком так, что
все области, на которые разбит лист складками, остаются плоскими, а перегибание происходит только по складкам.

В такой постановке правильный ответ --- «можно».
Удивительнее то, что возможность немного увеличить периметр видна в построении одной 
из самых древних фигурок оригами --- традиционного японского журавлика%
\footnote{Первая дошедшая до нас книга по оригами «Сембацуру ориката» была издана в 1797 году. 
Вся эта книга посвящена только складыванию журавлика.
Без сомнения, эта фигурка оригами гораздо древнее этой книги.}.
Про то, что оригамисты знали ответ, математики узнали в 1997 году\footnote{Вот \textattachfile[color=0 0 1,mimetype=text/html]{pics/napkin.html}{этот html-файл} расскажет, как это произошло.}.

Журавлика складывают из специальной заготовки, 
которая имеет четыре больших конца и ещё один короткий. 
Из двух концов делаются два крыла, а два других «утоньшаются»,
из них получаются хвост и шея.
Если же эту «операцию утоньшения» применить дважды к каждому большому концу, 
то заготовку можно будет выложить на плоскости так, 
что периметр превысит периметр исходного квадрата.

На рисунке~\ref{primer} справа изображена полученная фигура.
Подобный пример описан в книге Лэнга\footnote{R. Lang, \textit{Origami Design Secrets}}.
Если обозначить через $a$ сторону исходного квадрата, 
то обход каждой из четырёх игл займёт чуть меньше $a$, 
а обход короткого тупого конца займёт около $(\sqrt2-1){\cdot}a$, 
именно за счёт этого последнего конца мы получаем прибавку к периметру. 
Если много раз повторить процедуру утоньшения, то можно добиться, 
чтобы эта добавка стала произвольно близка к $(\sqrt2-1){\cdot}a\approx\tfrac25{\cdot}a$.

В полученной фигуре периметр увеличивается меньше чем на процент,
и сделать её руками непросто, хотя возможно.
Следует найти очень тонкую бумагу, так как максимальное количество слоёв после
расправления концов станет равно восьмидесяти.

\begin{figure}[h]
\ \ \ \ \ \begin{lpic}[t(0mm),b(0mm),r(0mm),l(0mm)]{pics/zagotovka(0.25)}
{\large
\lbl{230,170,-10;$\longrightarrow$}
\lbl{295,105,-20;$\longrightarrow$}
}
{\tiny
\lbl{119,119;$1$}
\lbl{134,110;$2$}
\lbl{147,110;$3$}
\lbl{160,110;$4$}
\lbl{170,110;$5$}
\lbl{180,110;$6$}
\lbl{190,110;$7$}
\lbl{200,110;$8$}

\lbl{117,100;$16$}
\lbl{134,100;$15$}
\lbl{147,100;$14$}
\lbl{160,100;$13$}
\lbl{170,100;$12$}
\lbl{180,100;$11$}
\lbl{190,100;$10$}
\lbl{200,100;$9$}

\lbl{110,134;$2$}
\lbl{110,147;$3$}
\lbl{110,160;$4$}
\lbl{110,171;$5$}
\lbl{110,181;$6$}
\lbl{110,191;$7$}
\lbl{110,201;$8$}

\lbl{98,120;$16$}
\lbl{98,134;$15$}
\lbl{98,147;$14$}
\lbl{98,160;$13$}
\lbl{98,171;$12$}
\lbl{98,181;$11$}
\lbl{98,191;$10$}
\lbl{98,201;$9$}

\lbl{90,110;$17$}
\lbl{74,110;$18$}
\lbl{60,110;$19$}
\lbl{48,110;$20$}
\lbl{37,110;$21$}
\lbl{26.5,110;$22$}
\lbl{17,110;$23$}
\lbl{7,110;$24$}

\lbl{74,100;$31$}
\lbl{60,100;$30$}
\lbl{48,100;$29$}
\lbl{37,100;$28$}
\lbl{27,100;$27$}
\lbl{17,100;$26$}
\lbl{7,100;$25$}

\lbl{110,88;$17$}
\lbl{110,74;$18$}
\lbl{110,59;$19$}
\lbl{110,48;$20$}
\lbl{110,38;$21$}
\lbl{110,28;$22$}
\lbl{110,18;$23$}
\lbl{110,8;$24$}

\lbl{88,88;$32$}

\lbl{98,74;$31$}
\lbl{98,59;$30$}
\lbl{98,48;$29$}
\lbl{98,38;$28$}
\lbl{98,28;$27$}
\lbl{98,18;$26$}
\lbl{98,8;$25$}
}
\end{lpic}
\caption{\label{primer} В центре --- сложенная заготовка журавлика с дважды заострёнными концами. 
Справа --- заготовка журавлика с расправленными концами.
Слева --- сетка складок; углы между всеми складками при вершинах квадрата равны; на каждом треугольнике указан номер его слоя в заготовке; складки, которые добавляются при расправлении концов, помечены точечным пунктиром.}
\end{figure}

\section{Ёжик и расчёска}

Оказывается, что периметр сложенной фигуры можно сделать произвольно большим.
Это можно увидеть из построения другой фигурки оригами, так называемого «морского ежа».
Эту фигурку сложил в 1987 году американский оригамист Роберт Лэнг.
Не зная ни про ёжика, ни про журавлика, Алексей Тарасов нашёл своё решение в 2004 году,%
\footnote{А. Тарасов, Решение задачи Арнольда о «мятом рубле», \textit{Чебышевский сборник} 5, № 1(9), 2004, 174---187.}
сложив другую фигуру, которую мы назовём расчёской Тарасова.
Он подошёл к задаче как математик, привёл точную формулировку и показал, что в процессе подобного складывания можно обойтись без растягиваний и сжиманий листа (приёмы, которыми оригамисты часто пользуются).

\begin{figure}[h]
\ \ \ \ 
\includegraphics[scale=0.28]{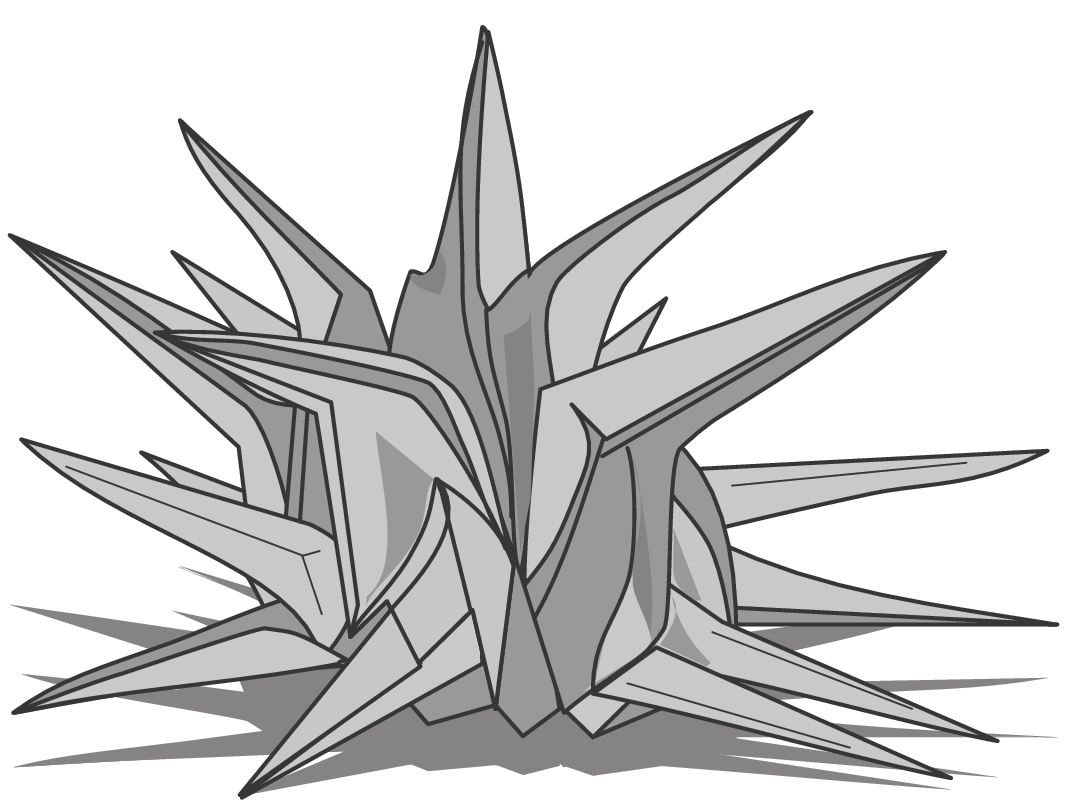}
\hfill
\includegraphics[scale=0.55]{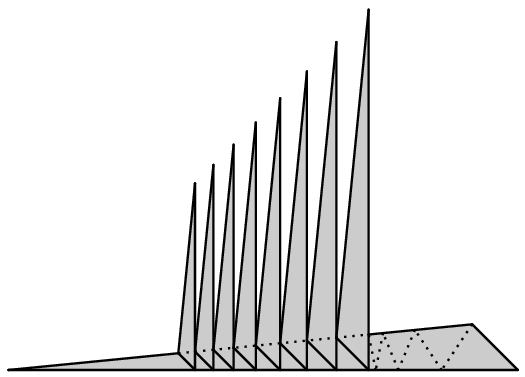}
\ \ \ \ 
\caption{\label{taras}Лэнговский ёж и расчёска Тарасова.}
\end{figure}

Эту «точную формулировку» я здесь не привожу, и она совсем не такая простая, как может показаться.
Объяснять процессы складывания я тоже не буду, их можно найти соответственно в книжке Монтролла и Лэнга%
\footnote{Montroll and Lang, \textit{Origami Sea Life}} 
и в видеофайлах Алексея Тарасова%
\footnote{Фильмы: складывание \textattachfile[color=0 0 1,author=Alexey Tarasov,mimetype=video/avi]{pics/rouble.avi}{расчёски}, отдельно \textattachfile[color=0 0 1,author=Alexey Tarasov,mimetype=video/avi]{pics/layer-up.avi}{одной иглы}
и \textattachfile[color=0 0 1,author=Alexey Tarasov,mimetype=video/avi]{pics/2layers.avi}{двух игл}.}. 
Но давайте поймём принцип, как можно сделать периметр произвольно большим.

Как видите, в обоих этих примерах имеется большое число иголок.
При складывании расчёски квадрат со стороной $a$ делится на $4{\cdot}n^2$ 
равных клеток --- $2{\cdot}n$ по вертикали и $2{\cdot}n$ по горизонтали 
(для расчёски на рисунке~\ref{taras} $n=2$).
Каждая из иголок получается из одной клетки, но используется только половина клеток; 
общее число игл равно $2{\cdot}n^2$. 
При этом толщину каждой иглы можно сделать произвольно малой, 
а длину больше, чем $\tfrac{a}{8{\cdot}n}$ --- четверти стороны клетки. 
Значит, обход каждой иглы займёт больше чем $\tfrac{a}{4{\cdot}n}$, 
и  общий периметр превысит $\tfrac{n{\cdot}a}2$.
Таким образом, взяв $n$ достаточно большим, 
можно добиться того, что периметр превзойдёт любое наперёд заданное число.

\medskip
При таком складывании размер расчёски будет меньше $\tfrac{a}{2{\cdot}n}$ --- стороны клетки, и простой подсчёт показывает, что наибольшее количество слоёв в одной точке должно быть существенно больше $n^4$ (иначе иглы будут находить друг на друга).
Поэтому этим способом практически невозможно существенно увеличить периметр обычного листа.
Представьте, например, что вы захотите таким способом увеличить в 10 раз периметр квадрата $10\times10$ сантиметров из обычного тетрадного листа.
Для этого придётся взять $n$ существенно больше $10$, размер полученной фигуры должен быть меньше сантиметра, и она будет более чем в 10 тысяч раз (т. е. $10^4$) толще листа.
Толщина тетрадного листа превышает 0,05 миллиметра\footnote{Толщину бумаги в тетради легко оценить, разделив её толщину на количество листов.}, а значит, наша фигура в сложенном виде должна быть «толще» половины метра!

\section{Ещё задачи}

Способы увеличения периметра, описанные выше, обладают одним недостатком --- при складывании необходимо одновременно сгибать лист в огромном числе мест.
Следующий пример основан на идее Ивана Ященко%
\footnote{см. Math. Intelligencer 1998 20(2) 38---40},
он также даёт лишь незначительное увеличение периметра, но состоит в последовательном применении одного и того же простого действия с добавлением пары складок за раз.
Если повторить его достаточно много раз, то получим фигуру, сколь угодно близкую к фигуре справа, у неё периметр слегка больше, чем у квадрата.

\begin{figure}[h]
\begin{lpic}[t(0mm),b(0mm),r(0mm),l(5mm)]{pics/yaschenko(0.57)}
\end{lpic}
\caption{\label{yasch} Сетка складок после седьмого повторения действия, сама складка и предельная фигура (после бесконечного числа повторений).}
\end{figure}

\begin{thm}{Задача}
Разберитесь в этом построении и докажите, что этот процесс действительно даёт увеличение периметра.

Оцените число повторений, после которых периметр становится больше периметра исходного квадрата. 
Найдите количество слоёв в самом толстом месте полученной фигуры.
\end{thm}

Предположим, квадрат сложен в виде плоской фигуры. 
Отметим все его складки и развернём. 
Складки дают разбиение квадрата на многоугольники, как, например, на первой картинке рисунка~\ref{yasch}.
Это разбиение называется \emph{сеткой складок}.

\begin{thm}{Задача}
Докажите, что 
многоугольники в сетке складок можно раскрасить в чёрный и белый цвета так, что соседние по стороне многоугольники имеют разный цвет,
и при этом сумма чёрных углов, сходящихся в одной внутренней вершине сетки, равна $180^\circ$.

Постройте разбиение квадрата на многоугольники такое, что для него выполняются оба указанных условия, но при этом оно не является сеткой складок.
\end{thm}

Расчёска Тарасова (как и фигуры в других подобных примерах) является невыпуклым многоугольником.
А что, если дополнительно потребовать, чтобы сложенная фигура была выпуклой? 
Это приводит к другой задаче с другим ответом.
Эта задача посложней, я привёл её полное решение в разделе~\ref{resh}.

\begin{thm}{Задача} 
\label{convex-prob}
Предположим, выпуклый многоугольник $M$ можно сложить на плоскости в виде другого выпуклого многоугольника $M'$. 
Докажите, что периметр $M'$ не превосходит периметра $M$.
(В частности, из квадрата нельзя сложить выпуклый многоугольник с б\'{о}льшим периметром).
\end{thm}

\medskip
А вот задача, которую мне решить не удалось.
Назовём фигуру $F$ \textit{вкладышем}, если она удовлетворяет следующему условию: 
\textit{любую плоскую фигуру, которую можно сложить из $F$, можно полностью накрыть копией~$F$}.
Другими словами, если фигура $F'$ сложена из $F$, то  $F'$ конгруэнтна подмножеству $F$.

\begin{wrapfigure}{r}{38mm}
\vskip-0mm
\centering
\includegraphics{mpost/pic-10}
\caption{}\label{counteexample}

\end{wrapfigure}

\begin{thm}{Задача}
Верно ли, что круги и только они являются вкладышами?%
\end{thm}

{\sloppy

\paragraph{P.S.} На рисунке \ref{counteexample} показан некруговой пример вкладыша, найденный позже.
Он ограничен дугой окружности и дугой параболы с общими концами, в которых они касаются друг друга.
Обоснование можно найти в обсуждении этой задачи на math\textit{overflow}.%
\footnote{\href{http://mathoverflow.net/questions/7016}{mathoverflow.net/questions/7016}.}
Очень похожий невыпуклый пример был найден раньше.%
\footnote{O. Aichholzer et al., \textit{Covering folded shapes},
J. Comput. Geom. 5 (2014), 150--168, Lemma~14.}

}
\par\WFclear

\section{Решение задачи 4}\label{resh}

Задача~\ref{convex-prob} также является переформулировкой частных случаев двух известных задач комбинаторной геометрии. 
Если не лень читать по-английски, посмотрите статьи Александера%
\footnote{R. Alexander, \textit{Lipschitzian mappings and total mean curvature of polyhedral surfaces I}, Trans. Amer. Math. Soc. 288 (1985), no. 2, 661--678.}
и Бездека --- Коннелли%
\footnote{K. Bezdek, R. Connelly, \textit{Pushing disks apart---the Kneser--Poulsen conjecture in the plane.}  J. Reine Angew. Math.  553  (2002), 221--236.}.
Там вы найдёте сами задачи и их решения,
оба этих решения очень красивые, но у меня приведено другое, более элементарное:

\medskip

\noi\textit{Решение.} Для точек $A$, $B$, $C$ в многоугольнике $M$ будем обозначать $A'$, $B'$, $C'$ соответствующие точки в многоугольнике $M'$.

Пусть $\{A_1,A_2,\dots,A_n\}$ --- конечное множество точек на плоскости.
Назовём обхватом $\{A_1,A_2,\dots,A_n\}$ наименьшую длину, которую надо пройти, чтобы обойти вокруг всех точек $A_1,A_2,\dots,A_n$.
Обхват $\{A_1,A_2,\dots,A_n\}$ будем обозначать $\ell(A_1,A_2,\dots,A_n)$.
Иначе говоря, $\ell(A_1,A_2,\dots,A_n)$ есть периметр наименьшего выпуклого многоугольника, содержащего все точки $A_1,A_2,\dots,A_n$.

\begin{wrapfigure}{r}{33mm}
\begin{lpic}[t(-2mm),b(0mm),r(0mm),l(0mm)]{pics/okhvat(0.6)}
\end{lpic}
\caption{Обхват точек}
\end{wrapfigure}

Будем рассуждать от противного,
предположим 
$$p(M')>p(M),$$
где $p(M)$ обозначает периметр $M$.
Тогда в $M$ можно найти множество из $n$ точек $\{A_1,A_2,\dots,A_n\}$ такое, 
что для соответствующего множества\linebreak
$\{A_1',A_2',\dots,A_n'\}$ в $M'$ выполняется неравенство
$$\ell(A_1',A_2',\dots,A_n')>\ell(A_1,A_2,\dots,A_n).$$
При этом можно предположить, что множество $\{A_1,A_2,\dots,A_n\}$ выбрано так, что
\begin{enumerate}
\item\label{min-n} число $n$ принимает минимально возможное значение;
\item\label{max-l} из всех $n$-точечных подмножеств $M$  величина
$$\ell(A_1',A_2',\dots,A_n')-\ell(A_1,A_2,\dots,A_n)>0$$
принимает максимально возможное значение.
\end{enumerate}
Давайте обозначим через $P$ и $P'$ наименьшие выпуклые многоугольники, 
содержащие множества $\{A_1,A_2,\dots,A_n\}$ и $\{A_1',A_2',\dots,A_n'\}$ соответственно.
В частности, 
$$p(P)=\ell(A_1,A_2,\dots,A_n)\ \ \ \text{и}\ \ \   p(P')=\ell(A_1',A_2',\dots,A_n').$$

Заметьте, что все точки $A_i'$ являются вершинами $P'$ и все $A_i'$ различные.
Действительно, предположим, $A_n'$ лежит внутри или на стороне $P'$ или \hbox{$A_n'=A_i'$} при $i\not=n$,
тогда 
$$\ell(A_1',A_2',\dots,A_{n-1}')= \ell(A_1',A_2',\dots,A_{n-1}',A_{n}').$$
При этом очевидно, что
$$\ell(A_1,A_2,\dots,A_{n-1})\le \ell(A_1,A_2,\dots,A_{n-1},A_{n}).$$
То есть мы получаем противоречие с условием \ref{min-n}.

Теперь заметьте, что если $A_i$ --- вершина $P$, то угол $P$ при $A_i$ 
не меньше угла $P'$ при $A_i'$. 
Действительно, обозначим угол $P$ при $A_i$ через $\alpha$, а 
угол $P'$ при $A_i'$ через $\alpha'$.
Если двигать $A_i$ с единичной скоростью внутрь $P$ вдоль биссектрисы угла, 
то обхват $\ell(A_1,A_2,\dots,A_{n})$
в начальный момент будет уменьшаться со скоростью $2{\cdot}\cos\tfrac\alpha2$.
При этом точка $A_i'$ будет двигаться в $M'$ также с единичной скоростью, 
и нетрудно видеть,
что в начальный момент обхват $\ell(A_1',A_2',\dots,A_{n}')$
не может уменьшаться со скоростью большей, чем $2{\cdot}\cos\tfrac{\alpha'}2$.
По условию \ref{max-l}, при таком движении разница
$$\ell(A_1',A_2',\dots,A_{n}')-\ell(A_1,A_2,\dots,A_{n})$$
не может увеличиваться.
Значит 
$$2{\cdot}\cos\tfrac{\alpha'}2
\ge
2{\cdot}\cos\tfrac\alpha2\ \ \text{или}\ \ \alpha'\le\alpha.$$

Применив теорему о сумме углов многоугольника, получаем, что и все $A_i$ являются вершинами $P$.\footnote{Мы также получаем, что углы при соответствующих вершинах равны. Это свойство в нашем доказательстве не применяется.}

Если предположить, что точки $A_i$ пронумерованы в порядке обхода $P$, то получаем:
$$\ell(A_1,A_2,\dots,A_{n})=A_1A_2+A_2A_3+\cdots+A_nA_1.$$
Очевидно, что $A_i'A_j'\le A_iA_j$ для всех $i$ и $j$, а значит
$$A_1'A_2'+A_2'A_3'+\cdots+A_n'A_1'\le A_1A_2+A_2A_3+\cdots+A_nA_1.$$
Далее, для любого конечного множества точек плоскости $\{A_1',A_2',\dots,A_{n}'\}$ выполняется
$$\ell(A_1',A_2',\dots,A_{n}')\le A_1'A_2'+A_2'A_3'+\cdots+A_n'A_1'.$$
Иначе говоря, обхват вершин любой замкнутой ломаной не превосходит длины ломаной.
Докажите это сами. 

Таким образом,
$$\ell(A_1',A_2',\dots,A_{n}')\le\ell(A_1,A_2,\dots,A_{n})$$
--- противоречие.\qed
\end{document}